# Intrinsic Location Parameter of a Diffusion Process

R. W. R. Darling[1][2]

---

[1] P. O. Box 535, Annapolis Junction, Maryland 20701-0535, USA. E-mail: rwrd@afterlife.ncsc.mil

[2] Supported by the Air Force Office of Scientific Research.

**Dedicated to James Eells on his 72nd Birthday**

ABSTRACT: For nonlinear functions $f$ of a random vector $Y$, $E[f(Y)]$ and $f(E[Y])$ usually differ. Consequently the mathematical expectation of $Y$ is not intrinsic: when we change coordinate systems, it is not invariant. This article is about a fundamental and hitherto neglected property of random vectors of the form $Y = f(X(t))$, where $X(t)$ is the value at time $t$ of a diffusion process $X$: namely that there exists a measure of location, called the "intrinsic location parameter" (ILP), which coincides with mathematical expectation only in special cases, and which is invariant under change of coordinate systems. The construction uses martingales with respect to the intrinsic geometry of diffusion processes, and the heat flow of harmonic mappings. We compute formulas which could be useful to statisticians, engineers, and others who use diffusion process models; these have immediate application, discussed in a separate article, to the construction of an intrinsic nonlinear analog to the Kalman Filter. We present here a numerical simulation of a nonlinear SDE, showing how well the ILP formula tracks the mean of the SDE for a Euclidean geometry.

RESUME: Pour une fonction non linéaire $f$ d'un vecteur aléatoire, $E[f(Y)]$ et $f(E[Y])$ sont usuellement différents. Par conséquent, l'espérance mathématique de $Y$ n'est pas intrinsèque: quand nous changeons le système des coordonnées, elle n'est pas invariante. Cet article concerne une propriété fondamentale, négligée jusqu'à maintenant, des vecteurs aléatoires de la forme $Y = f(X(t))$, où $X(t)$ est la valeur au temps $t$ d'un processus de diffusion $X$: c'est à dire qu'il existe une mesure de position, nommée le "paramètre intrinsèque de centrage" (PIC), qui coincide avec l'espérance mathématique seulement dans des cas spécifiques, et qui est invariante par changement du système des coordonnées. La construction utilise des martingales en rapport avec la géometrie intrinsèque des processus de diffusion, et le flot de chaleur des applications harmoniques. Nous calculons des formules qui peuvent être utiles aux statisticiens, aux ingénieurs, et à toute autre personne qui utilise des modèles fondés sur des processus de diffusion; ces formules se mettent en service à la construction d'une analogue non linéaire intrinsèque du filtre de Kalman, discutée dans un autre article. Nous présentons ici une simulation numérique d'une EDS non linéaire, qui montre la précision avec laquelle la formule de PIC suit la moyenne de l'EDS pour une géometrie Euclidéenne.

September 5, 1998 8:33 pm.

# 1    Introduction

## 1.1    Background

The relationship between martingales and parabolic partial differential equations was pointed out in the classic paper of Doob [11]: the solution $u(t, x)$ to the one-dimensional heat equation, with a given function $\psi$ as the boundary condition $u(0, .)$, is given by

$$u(T, x) = E_x[\psi(W_T)], \qquad (1)$$

where $\{W_t\}$ is Brownian motion with $W_0 \equiv x$. This can also be expressed as the initial value $V_0$ of the martingale $V_t \equiv u(T - t, W_t)$ which terminates at $\psi(W_T)$ at time $T$.

In the case of nonlinear parabolic PDE, the martingale $\{V_t\}$ must be replaced by the solution of an inverse problem for stochastic differential equations, also called a backwards SDE, as in the work of Pardoux and Peng: see [23], [10]. In the case of the system of elliptic PDE known as a harmonic mapping between Riemannian manifolds, this problem becomes one of constructing a martingale on a manifold with prescribed limit, which has been been solved in works by Kendall [20], [21], Picard [24], [25], Arnaudon [1], Darling [6], [7], and Thalmaier [29]. Thalmaier [30] studies the parabolic problem foir the nonlinear heat equation. The main point here is that the straightforward computation of an expectation as in (1) is no longer available in the nonlinear case. For discussion of the concept of using martingales on a manifold for determining barycentres, see Emery and Mokobodzki [18].

The aim of the present paper is to show the ideas mentioned in the previous paragraph have application to the question of determining the "mean" of a diffusion process, or of its image under a smooth function, in an intrinsic way, and that furthermore it is possible to compute an approximation to such a mean without excessive effort.

## 1.2    Main Results

Suppose $X$ is a Markov diffusion process on $R^p$, or more generally on a manifold $N$. The diffusion variance of $X$ induces a semi-definite metric $\langle .|. \rangle$ on the cotangent bundle, a version of the Levi-Civita connection $\Gamma$, and a Laplace-Beltrami operator $\Delta$. We may treat $X$ as a diffusion on $N$ with generator $\xi + (1/2)\Delta$, where $\xi$ is a vector field.

For sufficiently small $\delta > 0$, $X_\delta$ has an "intrinsic location parameter", defined to be the non-random initial value $V_0$ of a $\Gamma$-martingale $V$ terminating at $X_\delta$. It is obtained by solving a system of forward-backwards stochastic differential equations (FBSDE): a forward equation for $X$, and a backwards equation for $V$. This FBSDE is the stochastic equivalent of the heat equation (with drift $\xi$) for harmonic mappings, a well-known system of quasilinear PDE.

Let $\{\phi_t: N \to N, t \geq 0\}$ be the flow of the vector field $\xi$, and let $x_t \equiv \phi_t(x_0) \in N$. Our main result is that $\exp_{x_\delta}^{-1} V_0$ can be intrinsically approximated to first order in $T_{x_\delta} N$ by

$$\nabla d\phi_\delta(x_0)\,(\Pi_\delta) - \int_0^\delta (\phi_{\delta-t})_*(\nabla d\phi_t(x_0))\,d\Pi_t$$

where $\Pi_t = \int_0^t (\phi_{-s})_*\langle.|.\rangle_{x_s}\,ds \in T_{x_0}N \otimes T_{x_0}N$. This is computed in local coordinates. More generally, we find an intrinsic location parameter for $\psi(X_\delta)$, if $\psi: N \to M$ is a $C^2$ map into a Riemannian manifold $M$. We also treat the case where $X_0$ is random.

## 2   Geometry Induced by a Diffusion Process

### 2.1   Diffusion Process Model

Consider a Markov diffusion process $\{X_t, t \geq 0\}$ with values in a connected manifold $N$ of dimension $p$, represented in coordinates by

$$dX_t^i = b^i(X_t)\,dt + \sum_{j=1}^p \sigma_j^i(X_t)\,dW_t^j, \qquad (2)$$

where $\sum b^i \frac{\partial}{\partial x_i}$ is a vector field on $N$, $\sigma(x) \equiv (\sigma_j^i(x)) \in L(R^p; T_xN)$, and $W$ is a Wiener process in $R^p$. We assume for simplicity that the coefficients $b^i$, $\sigma_j^i$ are $C^2$ with bounded first derivative.

### 2.2   The diffusion variance semi-definite metric

Given a stochastic differential equation of the form (2) in each chart, it is well known that one may define a $C^2$ semi-definite metric $\langle.|.\rangle$ on the cotangent bundle, which we call the diffusion variance semi-definite metric, by the formula

$$\langle dx^i | dx^k \rangle_x \equiv (\sigma \cdot \sigma)^{ik}(x) \equiv \sum_{j=1}^p \sigma_j^i(x)\,\sigma_j^k(x). \qquad (3)$$

Note that $\langle.|.\rangle$ may be degenerate. This semi-definite metric is actually intrinsic: changing coordinates for the diffusion will give a different matrix $(\sigma_j^i)$, but the same semi-definite metric. We postulate:

*Axiom A:*   *The appropriate metric for the study of X is the diffusion variance semi-definite metric, not the Euclidean metric.*

The $p \times p$ matrix $((\sigma \cdot \sigma)^{ij})$ defined above induces a linear transformation $\alpha(x): T_x^*N \to T_xN$, namely

$$\alpha(x)(dx^i) \equiv \sum (\sigma \cdot \sigma)^{ij} \partial/\partial x_j.$$

Let us make a constant-rank assumption, i.e. that there exists a rank $r$ vector bundle $E \to N$, a subbundle of the tangent bundle, such that $E_x = \text{range}(\sigma(x)) \subseteq T_xN$ for all $x \in N$. In Section 7 below, we present a global geometric construction of what we call the canonical sub-Riemannian connection $\nabla^\circ$ for $\langle.|.\rangle$, with respect to a generalized inverse $g$, i.e. a vector bundle isomorphism $g: TN \to T^*N$ such that

$$\alpha(x) \bullet g(x) \bullet \alpha(x) = \alpha(x). \tag{4}$$

In local coördinates, $g(x)$ is expressed by a Riemannian metric tensor $(g_{rs})$, such that if $\alpha^{ij} \equiv (\sigma \cdot \sigma)^{ij}$, then

$$\sum_{r,s} \alpha^{ir} g_{rs} \alpha^{sj} = \alpha^{ij}. \tag{5}$$

The Christoffel symbols $\{\Gamma_{ij}^s\}$ for the canonical sub-Riemannian connection are specified by (84) below. The corresponding local connector $\Gamma(x) \in L(T_x R^p \otimes T_x R^p; T_x R^p)$ can be written in the more compact notation:

$$2g(\Gamma(x)(u \otimes v)) \cdot w = D\langle g(v)|g(w)\rangle(u) + D\langle g(w)|g(u)\rangle(v) - D\langle g(u)|g(v)\rangle(w), \tag{6}$$

where $g(\Gamma(x)(u \otimes v))$ is a 1-form, acting on the tangent vector $w$.

## 2.3   Intrinsic Description of the Process

The intrinsic version of (2) is to describe $X$ as a diffusion process on the manifold $N$ with generator

$$L \equiv \xi + \frac{1}{2}\Delta \tag{7}$$

where $\Delta$ is the (possibly degenerate) Laplace-Beltrami operator associated with the diffusion variance, and $\xi$ is a vector field, whose expressions in the local coordinate system $\{x^1, ..., x^p\}$ are as follows:

$$\Delta = \sum_{i,j} (\sigma \cdot \sigma)^{ij} \{D_{ij} - \sum_k \Gamma_{ij}^k D_k\}, \quad \xi = \sum_k \{b^k + \frac{1}{2}\sum_{i,j} (\sigma \cdot \sigma)^{ij} \Gamma_{ij}^k\} D_k. \tag{8}$$

Note that $\sum (\sigma \cdot \sigma)^{ij} \Gamma_{ij}^k$ has been specified by (3) and (6).

## 3    Γ-Martingales

Let $\bar{\Gamma}$ be a connection on a manifold $M$. An $H^2$ Γ-martingale is a kind of continuous semimartingale on $M$ which generalizes the notion of continuous $L^2$ martingale on $R^q$: see Emery [17] and Darling [7]. We summarize the main ideas, using global coordinates for simplicity.

Among continuous semimartingales in $R^q$, Itô's formula shows that local martingales are characterized by

$$f(X_t) - f(X_0) - (1/2)\int_0^t D^2 f(X_s)(dX \otimes dX)_s \in M_{loc}^c, \quad \forall f \in C^2(R^q), \tag{9}$$

where $(dX \otimes dX)^{ij}$ is the differential of the joint quadratic variation process of $X^i$ and $X^j$, and $M_{loc}^c$ refers to the space of real-valued, continuous local martingales (see Revuz and Yor [26]). For vector fields $\xi, \zeta$ on $R^q$, and $\omega \in \Omega^1(R^q)$, the smooth one-forms, a connection Γ gives an intrinsic way of differentiating $\omega$ along $\xi$ to obtain

$$\nabla_\xi \omega \in \Omega^1(R^q) \,.$$

$\nabla_\xi \omega \cdot \zeta$ is also written $\nabla \omega (\xi \otimes \zeta)$. When $\omega = df$, this gives the Hessian

$$\nabla df (D_k \otimes D_i) = D_{ki} f - \sum_j \Gamma^j_{ki} D_j f$$

where the $\{\Gamma^i_{jk}\}$ are the Christoffel symbols. The intrinsic, geometric restatement of (9) is to characterize a $\Gamma$–martingale $X$ by the requirement that

$$f(X_t) - f(X_0) - (1/2) \int_0^t \nabla df(X_s) (dX \otimes dX)_s \in M^c_{loc}, \; \forall f \in C^2(R^q) \,. \tag{10}$$

This is equivalent to saying that $M^k \in M^c_{loc}$ for $k = 1, \ldots, p$, where

$$dM^k_t = dX^k_t + (1/2) \sum_{i,j} \Gamma^k_{ij}(X_t) d\langle X^i, X^j \rangle_t \,. \tag{11}$$

If $N$ has a metric $g$ with metric tensor $(g_{ij})$, we say that $X$ is an $H^2$ $\Gamma$-martingale if (10) holds and also

$$E \langle X, X \rangle_\infty \equiv E \int_0^\infty \sum_{i,j} g_{ij}(X_t) d\langle X^i, X^j \rangle_t < \infty \,. \tag{12}$$

The $\Gamma$-martingale Dirichlet problem, which has been studied by, among others, Emery [16], Kendall [20], [21], Picard [24], [25], Arnaudon [1], Darling [6], [7], and Thalmaier [29], [30], is to construct a $\Gamma$-martingale, adapted to a given filtration, and with a given terminal value; for the Euclidean connection this is achieved simply by taking conditional expectation with respect to every $\sigma$-field in the filtration, but for other connections this may be as difficult as solving a system of nonlinear partial differential equations, as we shall now see.

## 4    Taylor Approximation of a Gamma-Martingale Dirichlet Problem

### 4.1    Condition for the Intrinsic Location Parameter

Consider a diffusion process $\{X_t, 0 \leq t \leq \delta\}$ on a $p$-dimensional manifold $N$ with generator $\xi + \frac{1}{2}\Delta$ where $\Delta$ is the Laplace-Beltrami operator associated with the diffusion variance, and $\xi$ is a vector field, as in (8). The coordinate-free construction of the diffusion $X$, given a Wiener process $W$ on $R^p$, uses the linear or orthonormal frame bundle: see Elworthy [15] p. 252. We suppose $X_0 = x_0 \in N$.

Also suppose $(M, h)$ is a Riemannian manifold, with Levi-Civita connection $\bar{\Gamma}$, and $\psi: N \to M$ is a $C^2$ map. The case of particular interest is when $M = N$, $\psi = \text{identity}$, and the metric on $N$ is a "generalized inverse" to $\sigma \cdot \sigma$ in the sense of (5). The general case of $\psi: N \to M$ is needed in the context of nonlinear filtering: see Darling [8].

Following Emery and Mokobodzki [18], we assert the following:

**Axiom B:** *Any intrinsic location parameter for $\psi(X_\delta)$ should be the initial value $V_0$ of an $\{\mathfrak{I}_t^W\}$-adapted $H^2$ $\bar{\Gamma}$-martingale $\{V_t, 0 \leq t \leq \delta\}$ on M, with terminal value $V_\delta = \psi(X_\delta)$.*

This need not be unique, but we will specify a particular choice below. In the case where $\sigma \cdot \sigma(x)$ does not depend on $x$, then the local connector $\Gamma$, given by (6), is zero, and $V_0$ is simply $E[\psi(X_\delta)]$. However our assertion is that, when $\Gamma$ is not the Euclidean connection, the right measure of location is $V_0$, and not $E[\psi(X_\delta)]$. We begin by indicating why an exact determination of $V_0$ is not computationally feasible in general.

## 4.2  Relationship with Harmonic Mappings

For simplicity of exposition, let us assume that there are diffeomorphisms $\varphi: N \to R^p$ and $\bar{\varphi}: M \to R^q$ which induce global coordinate systems $\{x^1, ..., x^p\}$ for N and $\{y^1, ..., y^q\}$ for M, respectively. By abuse of notation, we will usually neglect the distinction between $x \in N$ and $\varphi(x) \in R^p$, and write $x$ for both. $\Gamma(x)((\sigma \cdot \sigma)(x)) \in T_x R^p$ is given by (3) and (6), and the local connector $\bar{\Gamma}(y) \in L(T_y R^q \otimes T_y R^q; T_y R^q)$ comes from the Levi-Civita connection for $(M, h)$.

In order to find $\{V_t\}$, we need to construct an auxiliary adapted process $\{Z_t\}$, with values in $L\left(R^p; T_{V_t} R^q\right)$, such that the processes $\{X_t\}$ and $\{(V_t, Z_t)\}$ satisfy the following system of forward-backwards SDE:

$$X_t = x_0 + \int_0^t b(X_s)\, ds + \int_0^t \sigma(X_s)\, dW_s, \quad 0 \leq t \leq \delta; \tag{13}$$

$$V_t = \psi(X_\delta) - \int_t^\delta Z_s\, dW_s + \frac{1}{2}\int_t^\delta \bar{\Gamma}(V_s)(Z_s \cdot Z_s)\, ds, \quad 0 \leq t \leq \delta. \tag{14}$$

We also require that

$$E\left[\int_0^\delta \sum_{i,j} h_{ij}(V_s)(Z_s \cdot Z_s)^{ij} ds\right] < \infty. \tag{15}$$

[Equation (14) and condition (15) together say that $V$ is an $H^2$ $\bar{\Gamma}$-martingale, in the sense of (11) and (12).] Such systems are treated by Pardoux and Peng [23], but existence and uniqueness of solutions to (14) are outside the scope of their theory, because the coefficient $\bar{\Gamma}(v)(z \cdot z)$ is not Lipschitz in $z$.

However consider the second fundamental form $\nabla d\phi$ of a $C^2$ mapping $\phi: N \to M$. Recall that $\nabla d\phi(x) \in L(T_x N \otimes T_x N; T_{\phi(x)} M)$ may be expressed in local coordinates by:

$$\nabla d\phi(x)(v \otimes w) = D^2\phi(x)(v \otimes w) - D\phi(x)\Gamma(x)(v \otimes w) + \bar{\Gamma}(y)(D\phi(x)v \otimes D\phi(x)w) \tag{16}$$

for $(v, w) \in T_x R^p \times T_x R^p$, $y \equiv \phi(x)$. Let $\xi$ be as in (8). Consider a system of quasilinear parabolic PDE (a "heat equation with drift" for harmonic mappings - see Eells and Lemaire [13], [14]) consisting of a suitably differentiable family of mappings $\{u(t, .): N \to M\}$, for $t \in [0, \delta]$, such that

$$\frac{\partial u}{\partial t} = du \cdot \xi + \frac{1}{2}\nabla du\,(\sigma \cdot \sigma)\,,\ 0 \le t \le \delta, \tag{17}$$

$$u\,(0,\,.) = \psi. \tag{18}$$

For $x \in N$, the right side of (17) is $du\,(t,\,.) \cdot \xi(x) + \frac{1}{2}\nabla du\,(t,\,.)\,(\sigma \cdot \sigma(x)) \in T_{u(t,x)}M$. Following the approach of Pardoux and Peng [23], Itô's formula shows that

$$V_t = u\,(\delta - t, X_t) \in M, \tag{19}$$

$$Z(t) = du\,(\delta - t, X_t) \bullet \sigma(X_t) \in L\!\left(R^p; T_{V_t}M\right) \tag{20}$$

solves (14). In particular $u\,(\delta, x_0) = V_0$. (A similar idea was used by Thalmaier [29].)

### 4.2.a Comments on the Local Solvability of (17) - (18)

Recall that the energy density of $\psi : N \to M$ is given by

$$e\,(\psi)\,(x) \equiv \frac{1}{2}\|d\psi \otimes d\psi\,(\sigma \cdot \sigma)\|^2_{\psi(x)} = \frac{1}{2}\sum_{\beta,\gamma} h_{\beta\gamma}(\psi(x))\,d\psi^{\beta}(x) \otimes d\psi^{\gamma}(x)\,(\sigma \cdot \sigma(x))\,. \tag{21}$$

Note, incidentally, that this formula still makes sense when $\sigma \cdot \sigma$ is degenerate. In the case where $\sigma \cdot \sigma$ is non-degenerate and smooth, $\xi = 0$, and $\int e\,(\psi)\,d\,(\mathrm{vol}_N) < \infty$, the inverse function theorem method of Hamilton [19], page 122, suffices to show existence of a unique smooth solution to (17) - (18) when $\delta > 0$ is sufficiently small. For a more detailed account of the properties of the solutions when $\dim(N) = 2$, see Struwe [28], pages 221 - 235. Whereas Eells and Sampson [12] showed the existence of a unique global solution when $(M, h)$ has non-positive curvature, Chen and Ding [4] showed that in certain other cases blow-up of solutions is inevitable. The case where $\sigma \cdot \sigma$ is degenerate appears not to have been studied in the literature of variational calculus, and indeed is not within the scope of the classical PDE theory of Ladyzenskaja, Solonnikov, and Ural'ceva [22]. A probabilistic construction of a solution, which may or may not generalize to the case where $\sigma \cdot \sigma$ is degenerate, will appear in Thalmaier [30]. Work by other authors, using Hörmander conditions on the generator $\xi + \frac{1}{2}\Delta$, is in progress. For now we shall merely assume:

***Hypothesis I***    *Assume conditions on $\xi$, $\sigma \cdot \sigma$, $\psi$, and $h$ sufficient to ensure existence and uniqueness of a solution $\{u\,(t,\,.) : N \to M, 0 \le t \le \delta_1\}$, for some $\delta_1 > 0$.*

### 4.3    Definition: the Intrinsic Location Parameter

*For $0 \le \delta \le \delta_1$, the intrinsic location parameter of $\psi(X_\delta)$ is defined to be $u\,(\delta, x_0)$, where $x_0 = X_0$.*

This depends upon the generator $\xi + \frac{1}{2}\Delta$, given in (8), where $\Delta$ may be degenerate; on the mapping $\psi : N \to M$; and on the metric $h$ for $M$. It is precisely the initial value of an $\{\mathfrak{I}_t^W\}$ -adapted $H^2$ $\bar{\Gamma}$-martingale on M, with terminal value $V_\delta = \psi(X_\delta)$. However by using the solution of the PDE, we force the intrinsic location parameter to be unique, and to have some regularity as a function of $x_0$.

The difficulty with Definition 4.3 is that, in filtering applications, it is not feasible to compute solutions to (17) and (18) in real time. Instead we compute an approximation, as we now describe.

### 4.4  A Parametrized Family of Heat Flows

Consider a parametrized family $\{u^\gamma, 0 \leq \gamma \leq 1\}$ of equations of the type (17), namely

$$\frac{\partial u^\gamma}{\partial t} = du^\gamma \cdot \xi + \frac{\gamma}{2} \nabla du^\gamma (\sigma \cdot \sigma), \ 0 \leq t \leq \delta, \tag{22}$$

$$u^\gamma(0, .) = \psi. \tag{23}$$

Note that the case $\gamma = 1$ gives the system (17), while the case $\gamma = 0$ gives $u^0(t, x) = \psi(\phi_t(x))$, where $\{\phi_t, t \geq 0\}$ is the flow of the vector field $\xi$.

In a personal communication, Etienne Pardoux has indicated the possibility of a probabilistic construction, involving the system of FBSDE (47) and (56), of a unique family of solutions $\{u^\gamma\}$ for sufficiently small $\gamma \geq 0$, and for small time $\delta > 0$, based on the results of Darling [6] and methods of Pardoux and Peng [23]. For now, it will suffice to replace Hypothesis I by the following:

**Hypothesis II**   *Assume conditions on $\xi$, $\sigma \cdot \sigma$, $\psi$, and h sufficient to ensure existence of $\delta_1 > 0$ and $\gamma_1 > 0$ such that there is a unique $C^2$ mapping $(\gamma, t, x) \to u^\gamma(t, x)$ from $[0, \gamma_1] \times [0, \delta_1] \times N$ to M satisfying (22) and (23) for each $\gamma \in [0, \gamma_1]$.*

#### 4.4.a  Notation

For any vector field $\zeta$ on N, and any differentiable map $\phi: N \to P$ into a manifold P, the "push-forward" $\phi_*\zeta$ takes the value $d\phi \cdot \zeta(x) \in T_y P$ at $y \equiv \phi(x) \in P$; likewise $\phi_*(\zeta \otimes \zeta') \equiv \phi_*\zeta \otimes \phi_*\zeta'$.

We must also assume for the following theorem that we have chosen a generalized inverse $g: T^*N \to TN$ to $\sigma \cdot \sigma$, in the sense of (4), so that we may construct a canonical sub-Riemannian connection $\nabla^\circ$ for $\langle .|. \rangle$, with respect to g.

We now state the first result, which will later be subsumed by Theorem 4.7.

### 4.5  Theorem (PDE Version)

*Assume Hypothesis II, and that $0 \leq \delta \leq \delta_1$. Then, in the tangent space $T_{\psi(x_\delta)}M$,*

$$\left.\frac{\partial}{\partial \gamma}u^\gamma(\delta, x_0)\right|_{\gamma = 0} = \frac{1}{2}\{\nabla d\psi(x_\delta)(\phi_{\delta*}\Pi_\delta) + \psi_*\{\nabla^\circ d\phi_\delta(x_0)(\Pi_\delta) - \int_0^\delta (\phi_{\delta - t})_* \nabla^\circ d\phi_t(x_0) d\Pi_t\}\} \tag{24}$$

*where $\{\phi_t, t \geq 0\}$ is the flow of the vector field $\xi$, $x_t \equiv \phi_t(x_0)$, and*

$$\Pi_t \equiv \int_0^t (\phi_{-s})_* \langle .|. \rangle_{x_s} ds \in T_{x_0}N \otimes T_{x_0}N. \tag{25}$$

In the special case where $M = N$, $\psi$ = identity, and $h = g$, the right side of (24) simplifies to the part in parentheses {…}.

### 4.5.a    Definition

*The expression (24) is called the approximate intrinsic location parameter in the tangent space $T_{\psi(x_\delta)}M$, denoted $I_{x_0}[\psi(X_\delta)]$.*

### 4.5.b    Remark: How the Formula is Useful

First we solve the ODE for the flow $\{\phi_t, 0 \leq t \leq \delta\}$ of the vector field $\xi$, compute $\partial u^\gamma(\delta, x_0)/\partial \gamma$ at $\gamma = 0$ using (24) (or rather, using the local coordinate version (33)), then use the exponential map to project the approximate location parameter on to $N$, giving

$$\exp_{\psi(x_\delta)} \left\{ \frac{\partial}{\partial \gamma} u^\gamma(\delta, x_0) \bigg|_{\gamma = 0} \right\} \in M. \tag{26}$$

Computation of the exponential map likewise involves solving an ODE, namely the geodesic flow on $M$. In brief, we have replaced the task of solving a system of PDE by the much lighter task of solving two ODE's and performing an integration.

## 4.6    The Stochastic Version

We now prepare an alternative version of the Theorem, in terms of FBSDE, in which we give a local coordinate expression for the right side of (24). In this context it is natural to define a new parameter $\varepsilon$, so that $\gamma = \varepsilon^2$ in (22). Instead of $X$ in (13), we consider a family of diffusion processes $\{X^\varepsilon, \varepsilon \geq 0\}$ on the time interval $[0, \delta]$, where $X^\varepsilon$ has generator $\xi + \varepsilon^2 \Delta/2$. Likewise $V$ in (14) will be replaced by a family $\{V^\varepsilon, \varepsilon \geq 0\}$ of $H^2$ $\bar{\Gamma}$-martingales, with $V^\varepsilon_\delta = \psi(X^\varepsilon_\delta)$, and $V^{\sqrt{\gamma}}_0 = u^\gamma(\delta, x_0)$. Note, incidentally, that such parametrized families of $\bar{\Gamma}$-martingales are also treated in recent work of Arnaudon and Thalmaier [2], [3].

### 4.6.a    Generalization to the Case of Random Initial Value

Suppose that, instead of $X_0 = x_0 \in N$ as in (13), we have $X_0 = \exp_{x_0}(U_0)$, where $U_0$ is a zero-mean random variable in $T_{x_0}N$, independent of $W$, with covariance $\Sigma_0 \in T_{x_0}N \otimes T_{x_0}N$; the last expression means that, for any pair of cotangent vectors $\beta, \lambda \in T^*_{x_0}N$,

$$E[(\beta \cdot U_0)(\lambda \cdot U_0)] = (\beta \otimes \lambda) \cdot \Sigma_0.$$

Now set up the family of diffusion processes $\{X^\varepsilon, \varepsilon \geq 0\}$ with initial values

$$X^\varepsilon_0 = \exp_{x_0}(\varepsilon U_0). \tag{27}$$

Each $\bar{\Gamma}$-martingale $\{V^\varepsilon_t, 0 \leq t \leq \delta\}$ is now adapted to the larger filtration $\{\tilde{\mathfrak{I}}^W_t\} \equiv \{\mathfrak{I}^W_t \vee \sigma(U_0)\}$. In particular,

$$\exp^{-1}_{\psi(x_\delta)} V^\varepsilon_0$$

is now a random variable in $T_{\psi(x_\delta)}M$ depending on $U_0$.

### 4.6.b   Definition

*In the case of a random initial value $X_0$ as above, the approximate intrinsic location parameter of $\psi(X_\delta)$ in the tangent space $T_{\psi(x_\delta)} M$, denoted $I_{x_0, \Sigma_0}[\psi(X_\delta)]$, is defined to be*

$$\frac{\partial}{\partial(\varepsilon^2)} E\left[\exp^{-1}_{\psi(x_\delta)} V_0^\varepsilon\right]\bigg|_{\varepsilon=0}. \tag{28}$$

We will see in Section 6.3 below that this definition makes sense. This is the same as

$$\frac{\partial}{\partial \gamma} E\left[\exp^{-1}_{\psi(x_\delta)} u^\gamma\left(\delta, X_0^{\sqrt{\gamma}}\right)\right]\bigg|_{\gamma=0},$$

and coincides with $I_{x_0}[\psi(X_\delta)]$, given by (24), in the case where $\Sigma_0 = 0$.

### 4.6.c   Some Integral Formulas

Given the flow $\{\phi_t, 0 \le t \le \delta\}$ of the vector field $\xi$, the derivative flow is given locally by

$$\tau_s^t \equiv d(\phi_t \bullet \phi_s^{-1})(x_s) \in L(T_{x_s} N; T_{x_t} N), \tag{29}$$

where $x_s = \phi_s(x)$, for $0 \le s \le \delta$. In local coordinates, we compute $\tau_s^t$ as a $p \times p$ matrix, given by

$$\tau_s^t = \exp\left\{\int_s^t D\xi(x_u)\, du\right\}.$$

Introduce the deterministic functions

$$\chi_t \equiv \int_0^t (\phi_{t-s})_* \langle \cdot | \cdot \rangle_{x_s} ds = \int_0^t \tau_s^t (\sigma \cdot \sigma)(x_s)(\tau_s^t)^T ds \in T_{x_t} N \otimes T_{x_t} N, \tag{30}$$

$$\Xi_t \equiv \chi_t + \tau_0^t \Sigma_0 (\tau_0^t)^T \in T_{x_t} N \otimes T_{x_t} N. \tag{31}$$

Note, incidentally, that $\Xi_t$ could be called the intrinsic variance parameter of $X_t$.

### 4.7   Theorem (Local Coördinates, Random Initial Value Version)

*Under the conditions of Theorem 4.5, with random initial value $X_0$ as in Section 4.6.a, the approximate intrinsic location parameter $I_{x_0, \Sigma_0}[\psi(X_\delta)]$ exists and is equal to the right side of (24), after redefining*

$$\Pi_t \equiv \Sigma_0 + \int_0^t (\phi_{-s})_* \langle \cdot | \cdot \rangle_{x_s} ds \in T_{x_0} N \otimes T_{x_0} N. \tag{32}$$

*In local coordinates, $2I_{x_0, \Sigma_0}[\psi(X_\delta)]$ is given by*

$$J\int_0^\delta \tau_t^\delta [D^2\xi(x_t)(\Xi_t) - \Gamma(x_t)(\sigma \cdot \sigma(x_t))]\, dt + D^2\psi(x_\delta)(\Xi_\delta) - J\tau_0^\delta \Gamma(x_0)(\Sigma_0) + \bar{\Gamma}(y_\delta)(J\Xi_\delta J^T) \tag{33}$$

*where $y_\delta \equiv \psi(x_\delta)$, $J \equiv D\psi(x_\delta)$, and $\Xi_t$ is given by (31).*

### 4.7.a   Remarks

- Theorem 4.7 subsumes Theorem 4.5, which corresponds to the case $\Sigma_0 = 0$.

- In the special case where $M = N$, $\psi =$ identity, and $h = g$, formula (33) reduces to:

$$I_{x_0, \Sigma_0}[X_\delta] = \frac{1}{2} \{ \int_0^\delta \tau_t^\delta [D^2\xi(x_t)(\Xi_t) - \Gamma(x_t)(\sigma \cdot \sigma(x_t))] \, dt - \tau_0^\delta \Gamma(x_0)(\Sigma_0) + \Gamma(x_\delta)(\Xi_\delta) \}. \quad (34)$$

- In the filtering context [8], formulas (33) and (34) are of crucial importance.

## 5   Example of Computing an Intrinsic Location Parameter

The following example shows that Theorem 4.7 leads to feasible and accurate calculations.

### 5.1   Target Tracking

In target tracking applications, it is convenient to model target acceleration as an Ornstein-Uhlenbeck process, with the constraint that acceleration must be perpendicular to velocity. Thus $(v, a) \in R^3 \times R^3$ must satisfy $v \cdot a = 0$, and the trajectory must lie within a set on which $\|v\|^2$ is constant. Therefore we may identify the state space $N$ with $TS^2 \subset R^6$, since the $v$-component lies on a sphere, and the $a$-component is perpendicular to $v$, and hence tangent to $S^2$.

Within a Cartesian frame, $X$ is a process in $R^6$ with components $V$ (velocity) and $A$ (acceleration), and the equations of motion take the nonlinear form:

$$\begin{bmatrix} dV \\ dA \end{bmatrix} = \begin{bmatrix} 0_{3\times 3} & I_3 \\ -\rho(X)I_3 & -\lambda P(V) \end{bmatrix} \begin{bmatrix} V \\ A \end{bmatrix} dt + \begin{bmatrix} 0 \\ \gamma P(V) \, dW(t) \end{bmatrix}. \quad (35)$$

Here the square matrix consists of four $3 \times 3$ matrices, $\lambda$ and $\gamma$ are constants, $W$ is a three-dimensional Wiener process, and if $x^T \equiv (v^T, a^T)$,

$$\rho(x) \equiv \|a\|^2 / \|v\|^2, \quad (36)$$

$$P(v) \equiv I - \frac{vv^T}{\|v\|^2} \in L(R^3; R^3), \quad (37)$$

Note that $P(v)$ is precisely the projection onto the orthogonal complement of $v$ in $R^3$, and $\rho(x)$ has been chosen so that $d(V \cdot A) = 0$.

### 5.2   Geometry of the State Space

The diffusion variance metric (3) is degenerate here; noting that $P^2 = P$, we find

$$\alpha \equiv \sigma \cdot \sigma \equiv \begin{bmatrix} 0_{3\times 3} & 0_{3\times 3} \\ 0_{3\times 3} & \gamma^2 P(v) \end{bmatrix}. \quad (38)$$

The rescaled Euclidean metric $g = \gamma^{-2} I_6$ on $R^6$ is a generalized inverse to $\alpha$ in the sense of (4), since $P^2 = P$. We break down a tangent vector $\zeta$ to $R^6$ into two 3-dimensional components $\zeta_v$ and $\zeta_a$. The constancy of $\|v\|^2$ implies that

$$DP(v)\eta_v = \frac{-1}{\|v\|^2}\{\eta_v v^T + v\eta_v^T\}. \tag{39}$$

Referring to formula (6) for the local connector $\Gamma(x)$,

$$D\langle g(\zeta)|g(\varsigma)\rangle(\eta) = \frac{-1}{\gamma^2\|v\|^2}\zeta_a^T(\eta_v v^T + v\eta_v^T)\varsigma_a, \ (\zeta,\varsigma,\eta) \in T_xN \times T_xN \times T_xN.$$

Taking first and second derivatives of the constraint $v \cdot a = 0$, we find that

$$\zeta_a^T v + \zeta_v^T a = 0, \ \eta_a^T \zeta_v + \eta_v^T \zeta_a = 0. \tag{40}$$

Using the last identity, we obtain from (6) the formula

$$\Gamma(x)(\zeta \otimes \varsigma) = \frac{S(\zeta \otimes \varsigma)v}{2\|v\|^2}, \ S(\zeta \otimes \varsigma) \equiv \begin{bmatrix} \zeta_a\varsigma_a^T + \varsigma_a\zeta_a^T \\ -\zeta_v\varsigma_a^T - \varsigma_v\zeta_a^T \end{bmatrix}. \tag{41}$$

In order to compute (33), note that, in particular,

$$\Gamma(x)(\sigma \cdot \sigma(x)) = \frac{\gamma^2}{\|v\|^2}\begin{bmatrix} P(v) \\ 0 \end{bmatrix}v = \begin{bmatrix} 0 \\ 0 \end{bmatrix}. \tag{42}$$

## 5.3  Derivatives of the Dynamical System

It follows from (8), (35), and (42) that the formula for the intrinsic vector field $\xi$ is:

$$\xi(x) = \begin{bmatrix} a \\ -\rho(x)v - \lambda P(v)a \end{bmatrix}. \tag{43}$$

Differentiate under the assumptions $\|v\|^2$ is constant and $v \cdot a = 0$, to obtain

$$D\xi(x)(\zeta) = \begin{bmatrix} 0 & I \\ \lambda Q - \rho & -\lambda P - 2Q \end{bmatrix}\begin{bmatrix} \zeta_v \\ \zeta_a \end{bmatrix}, \ Q \equiv \frac{va^T}{\|v\|^2}. \tag{44}$$

Differentiating (44) and using the identities (40),

$$D^2\xi(x)(\eta \otimes \zeta) = \frac{-2}{\|v\|^2}\begin{bmatrix} 0 \\ \eta_a^T\zeta_a v + (\zeta_v\eta_a^T + \eta_v\zeta_a^T)a \end{bmatrix}. \tag{45}$$

### 5.3.a   Constraints in the Tangent Space of the State Manifold

Let us write a symmetric tensor $\chi \in T_x N \otimes T_x N$ in $3 \times 3$ blocks as the matrix

$$\begin{bmatrix} \chi_{vv} & \chi_{va} \\ \chi_{av} & \chi_{aa} \end{bmatrix},$$

where $\chi_{av}^T = \chi_{va}$. Replacing $\eta_v \zeta_a^T$ by $\chi_{va}$, and $\eta_a^T \zeta_a$ by $\text{Tr}(\chi_{aa})$, etc., in (45), we find that

$$D^2 \xi(x)(\chi) = \frac{-2}{\|v\|^2} \begin{bmatrix} 0 \\ \text{Tr}(\chi_{aa}) v + (\chi_{va} + \chi_{av}) a \end{bmatrix}. \tag{46}$$

### 5.4   Ingredients of the Intrinsic Location Parameter Formula

Let $\psi: N \to R^6$ be the inclusion of the state space $N \cong TS^2$ into Euclidean $R^6$. Thus in formula (33), the local connector $\bar{\Gamma}(.)$ is zero on the target manifold, $J$ is the identity, and $D^2 \psi$ is zero. When $\Sigma_0$ is taken to be zero, the formula for the approximate intrinsic location parameter $m_\delta \equiv I_{x_0}[X_\delta]$ for $X_\delta$ becomes:

$$m_t = \int_0^t \tau_u^t H_u \, du, \quad H_t \equiv \frac{-1}{\|v_0\|^2} \begin{bmatrix} 0 \\ \text{Tr}(\chi_{aa}(t)) v_t + (\chi_{va}(t) + \chi_{av}(t)) a_t \end{bmatrix}$$

where $v_t$ and $a_t$ are the velocity and acceleration components of $x_t \equiv \phi_t(x)$, for $0 \le t \le \delta$, and $\tau_u^t$ and $\chi_t$ are given by (29) and (30). A straightforward integration scheme for calculating $(\tau_u^t, \chi_t, m_t)$ at the same time, using a discretization of $[0, \delta]$, is:

$$\tau_u^t \approx \exp\left\{ \frac{t-u}{2} [D\xi(x_u) + D\xi(x_t)] \right\},$$

$$\chi_t \approx \frac{t-u}{2} (\sigma \cdot \sigma)(x_t) + \tau_u^t \left[ \chi_u + \frac{t-u}{2} (\sigma \cdot \sigma)(x_u) \right] (\tau_u^t)^T,$$

$$m_t \approx \tau_u^t \left[ m_u + \frac{t-u}{2} H_u \right] + \frac{t-u}{2} H_t.$$

Since the local connector is zero on the target manifold, geodesics are simply straight lines, and $x_\delta + I_{x_0}[X_\delta]$ is a suitable estimate of the mean position of $X_\delta$.

## 5.5    Simulation Results

**FIGURE 1**   SIMULATIONS OF THE MEAN OF AN SDE, VERSUS ITS APPROXIMATE ILP

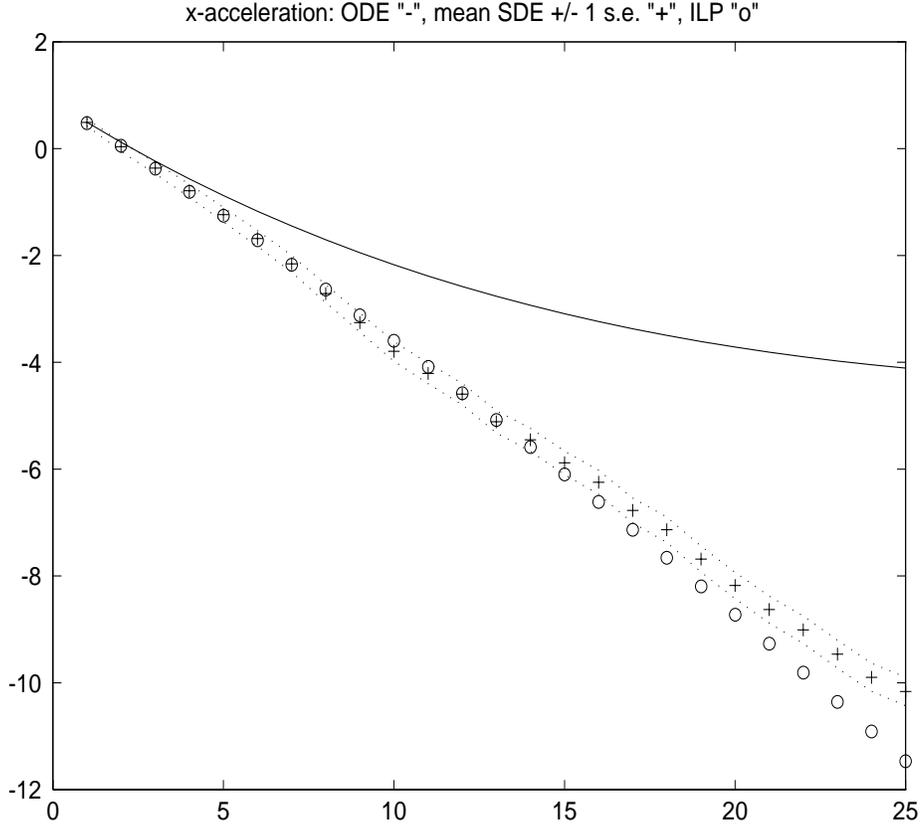

We created $10^4$ simulations of the process (35), with $\lambda = 0.5$ and $\gamma = 5.2 \times 10^3$, on the time interval $[0, 1]$, which was discretized into 25 subintervals for integration purposes. In each case $V$ and $A$ were initialized randomly, with magnitudes of 200 m/s and 50 m/s$^2$, respectively. The plot shows the $x$-component of acceleration (the other two are similar): the "+" signs represent the mean of the process (35) over $10^4$ simulations, with bands showing plus and minus one standard error, the solid line is the solution $\{x_t, 0 \leq t \leq 1\}$ of the discretized ODE $dx_t/dt = \xi(x_t)$, and the circles denote the approximate intrinsic location parameter (ILP). The reader will note that the ILP tracks the mean of the process better than the ODE does.

## 6    Proof of Theorem 4.7

The strategy of the proof will be to establish the formula (33) using Itô calculus, and then to show it is equivalent to (24) using differential geometric methods. While this may seem roundabout, the important formula for applications is really (33); converting it into (24) serves mainly as a check that formula (33) is indeed intrinsic. It will make no difference if we work in global coordinates, and identify $N$ with $R^p$ and $M$ with $R^q$.

## 6.1 Step I: Differentiation of the State Process with Respect to a Parameter

We consider a family of diffusion processes $\{X^\varepsilon, \varepsilon \geq 0\}$ on the time interval $[0, \delta]$, with initial values $X_0^\varepsilon = \exp_{x_0}(\varepsilon U_0)$; here $U_0$ is a zero-mean random variable in $T_{x_0} N$, independent of $W$, with covariance $\Sigma_0 \in T_{x_0} N \otimes T_{x_0} N$, and $X^\varepsilon$ has generator $\xi + \varepsilon^2 \Delta/2$.

Note that, in local coordinates, the SDE for $X^\varepsilon$ is not "$dX_s^\varepsilon = b(X_s^\varepsilon) ds + \varepsilon \sigma(X_s^\varepsilon) dW_s$", because the limiting case when $\varepsilon = 0$ would then be the ODE based on the vector field $\sum b^i D_i$, which is not the same as $\xi$, which is given by (8). Instead the SDE is

$$X_t^\varepsilon = \exp_{x_0}(\varepsilon U_0) + \int_0^t \xi(X_s^\varepsilon) ds - \varepsilon^2 \int_0^t \zeta(X_s^\varepsilon) ds + \int_0^t \varepsilon \sigma(X_s^\varepsilon) dW_s \qquad (47)$$

where we use the notation

$$\zeta(x) \equiv \frac{1}{2} \Gamma(x) (\sigma(x) \cdot \sigma(x)). \qquad (48)$$

In the case $\varepsilon = 0$, the solution is deterministic, namely $\{x_t, 0 \leq t \leq \delta\}$. Note that, in local coördinates,

$$\exp_{x_0}(\varepsilon U_0) = x_0 + \varepsilon U_0 - \frac{\varepsilon^2}{2} \Gamma(x_0)(U_0 \otimes U_0) + o(\varepsilon^2).$$

It is well known that, if the vector field $\xi$ and the semi-definite metric $\langle .|. \rangle$ are sufficiently differentiable, then the stochastic processes $\partial X^\varepsilon/\partial \varepsilon$ and $\partial^2 X^\varepsilon/\partial \varepsilon^2$ exist and satisfy the following SDEs:

$$\frac{\partial X_t^\varepsilon}{\partial \varepsilon} = U_0 + \int_0^t D\xi(X_s^\varepsilon) \left(\frac{\partial X_s^\varepsilon}{\partial \varepsilon}\right) ds + \int_0^t \sigma(X_s^\varepsilon) dW_s \qquad (49)$$

$$- 2\varepsilon \int_0^t \zeta(X_s^\varepsilon) ds + \varepsilon \int_0^t D\sigma(X_s^\varepsilon) \left(\frac{\partial X_s^\varepsilon}{\partial \varepsilon}\right) dW_s + O(\varepsilon^2);$$

$$\frac{\partial^2 X_t^\varepsilon}{\partial \varepsilon^2} = -\Gamma(x)(U_0 \otimes U_0) + \int_0^t D^2\xi(X_s^\varepsilon) \left(\frac{\partial X_s^\varepsilon}{\partial \varepsilon} \otimes \frac{\partial X_s^\varepsilon}{\partial \varepsilon}\right) ds + \int_0^t D\xi(X_s^\varepsilon) \left(\frac{\partial^2 X_s^\varepsilon}{\partial \varepsilon^2}\right) ds - 2\int_0^t \zeta(X_s^\varepsilon) ds \quad (50)$$

$$+ 2\int_0^t D\sigma(X_s^\varepsilon) \left(\frac{\partial X_s^\varepsilon}{\partial \varepsilon}\right) dW_s + O(\varepsilon),$$

where $O(\varepsilon)$ denotes terms of order $\varepsilon$. Define

$$\Lambda_t \equiv \left.\frac{\partial X_t^\varepsilon}{\partial \varepsilon}\right|_{\varepsilon = 0}, \quad \Lambda_t^{(2)} \equiv \left.\frac{\partial^2 X_t^\varepsilon}{\partial \varepsilon^2}\right|_{\varepsilon = 0}. \qquad (51)$$

Now (49) and (50) give:

$$d\Lambda_t = A_t \Lambda_t dt + \sigma(x_t) dW_t, \quad \Lambda_0 = U_0; \tag{52}$$

$$d\Lambda_t^{(2)} = [D^2\xi(x_t)(\Lambda_t \otimes \Lambda_t) + A_t \Lambda_t^{(2)} - 2\zeta(x_t)] dt + 2D\sigma(x_t)(\Lambda_t) dW_t, \tag{53}$$

$$\Lambda_0^{(2)} = -\Gamma(x_0)(U_0 \otimes U_0),$$

where $A_t \equiv D\xi(x_t)$. Let $\{\tau_s^t, 0 \leq s, t \leq \delta\}$ be the two-parameter semigroup of deterministic matrices given by (29), so that

$$\frac{\partial \tau_t^s}{\partial t} = -\tau_t^s A_t; \quad \tau_s^r = \tau_t^r \tau_s^t.$$

Then (52) becomes $d(\tau_t^0 \Lambda_t) = \tau_t^0 \sigma(x_t) dW_t$, which has a Gaussian solution

$$\Lambda_t = \tau_0^t U_0 + \tau_0^t \int_0^t \tau_s^0 \sigma(x_s) dW_s = \tau_0^t U_0 + \int_0^t \tau_s^t \sigma(x_s) dW_s. \tag{54}$$

Likewise (53) gives $d(\tau_t^0 \Lambda_t^{(2)}) = \tau_t^0 [D^2\xi(x_t)(\Lambda_t \otimes \Lambda_t) - 2\zeta(x_t)] dt + 2\tau_t^0 D\sigma(x_t)(\Lambda_t) dW_t$, whose solution is

$$\Lambda_t^{(2)} = -\tau_0^t \Gamma(x_0)(U_0 \otimes U_0) + \int_0^t \tau_s^t \{[D^2\xi(x_s)(\Lambda_s \otimes \Lambda_s) - 2\zeta(x_s)] ds + 2D\sigma(x_s)(\Lambda_s) dW_s\}. \tag{55}$$

## 6.2   Step II: Differentiation of the Gamma-Martingale with Respect to a Parameter

Consider the pair of processes $(V^\varepsilon, Z^\varepsilon)$ obtained from (19) and (20), where $u$ is replaced by $u^\varepsilon$. As in the case where $\varepsilon = 1$, $(V^\varepsilon, Z^\varepsilon)$ gives an adapted solution to the backwards equation corresponding to (14), namely

$$V_t^\varepsilon = \psi(X_\delta^\varepsilon) - \int_t^\delta Z_s^\varepsilon dW_s + \frac{1}{2} \int_t^\delta \bar{\Gamma}(V_s^\varepsilon)(Z_s^\varepsilon \cdot Z_s^\varepsilon) ds.$$

However the version of (20) which applies here is $Z_t^\varepsilon = Du^\varepsilon(\delta - t, X_t) \varepsilon \sigma(X_t)$, so we may replace $Z_s^\varepsilon$ by $\varepsilon Z_s^\varepsilon$, and the equation becomes

$$V_t^\varepsilon = \psi(X_\delta^\varepsilon) - \varepsilon \int_t^\delta Z_s^\varepsilon dW_s + \frac{\varepsilon^2}{2} \int_t^\delta \bar{\Gamma}(V_s^\varepsilon)(Z_s^\varepsilon \cdot Z_s^\varepsilon) ds. \tag{56}$$

By the regularity of $u^\varepsilon$, it follows that $\partial V^\varepsilon/\partial\varepsilon$, $\partial Z^\varepsilon/\partial\varepsilon$, $\partial^2 V^\varepsilon/\partial\varepsilon^2$, and $\partial^2 Z^\varepsilon/\partial\varepsilon^2$ exist, and satisfy the following equations:

$$\frac{\partial V_t^\varepsilon}{\partial \varepsilon} = D\psi(X_\delta^\varepsilon)\left(\frac{\partial X_\delta^\varepsilon}{\partial \varepsilon}\right) - \int_t^\delta Z_s^\varepsilon dW_s - \varepsilon \int_t^\delta \frac{\partial Z_s^\varepsilon}{\partial \varepsilon} dW_s + \varepsilon \int_t^\delta \bar{\Gamma}(V_s^\varepsilon)(Z_s^\varepsilon \cdot Z_s^\varepsilon) ds + O(\varepsilon^2); \tag{57}$$

$$\frac{\partial^2 V_t^\varepsilon}{\partial \varepsilon^2} = D\psi(X_\delta^\varepsilon)\left(\frac{\partial^2 X_\delta^\varepsilon}{\partial \varepsilon^2}\right) + D^2\psi(X_\delta^\varepsilon)\left(\frac{\partial X_\delta^\varepsilon}{\partial \varepsilon} \otimes \frac{\partial X_\delta^\varepsilon}{\partial \varepsilon}\right) - 2\int_t^\delta \frac{\partial Z_s^\varepsilon}{\partial \varepsilon} dW_s + \int_t^\delta \bar{\Gamma}(V_s^\varepsilon)(Z_s^\varepsilon \cdot Z_s^\varepsilon) ds + O(\varepsilon) \tag{58}$$

Note also that $V_t^0 = y_\delta \equiv \psi(x_\delta)$ for all $t \in [0, \delta]$. Take $\tilde{\mathfrak{F}}_t^W \equiv \mathfrak{F}_t^W \vee \sigma(U_0)$. By combining (54) and (57), we see that, if $J \equiv D\psi(x_\delta)$,

$$\Theta_t \equiv \left.\frac{\partial V_t^\varepsilon}{\partial \varepsilon}\right|_{\varepsilon=0} = E\left[J\Lambda_\delta \big| \tilde{\mathfrak{F}}_t^W\right] = J\tau_0^\delta U_0 + J\int_0^t \tau_s^\delta \sigma(x_s)\, dW_s; \tag{59}$$

$$Z_s^0 = J\tau_s^\delta \sigma(x_s). \tag{60}$$

Define

$$\Theta_t^{(2)} \equiv \left.\frac{\partial^2 V_t^\varepsilon}{\partial \varepsilon^2}\right|_{\varepsilon=0} \tag{61}$$

$$= E\left[J\Lambda_\delta^{(2)} + D^2\psi(x_\delta)(\Lambda_\delta \otimes \Lambda_\delta) + \int_t^\delta \bar{\Gamma}(y_\delta)(Z_s^0 \cdot Z_s^0)\, ds \,\Big|\, \tilde{\mathfrak{F}}_t^W\right] \tag{62}$$

From (55) and (62) we obtain:

$$\Theta_0^{(2)} = -J\tau_0^\delta \Gamma(x_0)(U_0 \otimes U_0) + JE\left[\left\{\int_0^\delta \tau_s^\delta [D^2\xi(x_s)(\Lambda_s \otimes \Lambda_s) - 2\zeta(x_s)]\, ds\right\} \Big| U_0\right]$$

$$+ E\left[\{D^2\psi(x_\delta)(\Lambda_\delta \otimes \Lambda_\delta) + \int_0^\delta \bar{\Gamma}(y_\delta)(Z_s^0 \cdot Z_s^0)\, ds\} \Big| U_0\right].$$

The expected value of a quadratic form $\eta^T A\eta$ in an $N_p(\mu, \Sigma)$ random vector $\eta$ is easily computed to be $\sum A_{ij} \Sigma^{ij} + \mu^T A\mu$. In this case,

$$\eta = \Lambda_s \sim N_p(\tau_0^t U_0, \chi_s),$$

where $\chi_t$ is given by (30), so we obtain

$$\Theta_0^{(2)} = -J\tau_0^\delta \Gamma(x_0)(U_0 \otimes U_0) + J\left[\int_0^\delta \tau_s^\delta [D^2\xi(x_s)(\chi_s) + D^2\xi(x_s)(\tau_0^s U_0 \otimes \tau_0^s U_0) - 2\zeta(x_s)]\, ds\right]$$

$$+ D^2\psi(x_\delta)(\chi_\delta) + D^2\psi(x_\delta)(\tau_0^\delta U_0 \otimes \tau_0^\delta U_0) + \bar{\Gamma}(y_\delta)(J\chi_\delta J^T). \tag{63}$$

## 6.3  Step III: A Taylor Expansion Using the Exponential Map

Let $y_\delta \equiv \psi(x_\delta)$, and define $\beta(\varepsilon) \equiv E\left[\exp_{y_\delta}^{-1} V_0^\varepsilon\right]$. Referring to (28), we are seeking $\left.\partial \beta/\partial \varepsilon^2\right|_{\varepsilon=0}$. It follows immediately from the geodesic equation that

$$D(\exp_y)^{-1}(y)(w) = w,\quad D^2(\exp_y)^{-1}(y)(v \otimes w) = \bar{\Gamma}(y)(v \otimes w) \tag{64}$$

It follows from (59) and (63) that

$$V_0^\varepsilon = y_\delta + \varepsilon J\tau_0^\delta U_0 + \frac{\varepsilon^2}{2}\Theta_0^{(2)} + o(\varepsilon^2) .$$

A Taylor expansion based on (64) gives

$$\exp_{y_\delta}^{-1} V_0^\varepsilon = V_0^\varepsilon - y_\delta + \frac{1}{2}\bar{\Gamma}(y_\delta)((V_0^\varepsilon - y_\delta) \otimes (V_0^\varepsilon - y_\delta)) + o\left(\left\|V_0^\varepsilon - y_\delta\right\|^2\right)$$

$$= \varepsilon J\tau_0^\delta U_0 + \frac{\varepsilon^2}{2}\{\Theta_0^{(2)} + \bar{\Gamma}(y_\delta)(J\tau_0^\delta U_0 \otimes J\tau_0^\delta U_0)\} + o(\varepsilon^2) .$$

Taking expectations, and recalling that $U_0$ has mean zero and covariance $\Sigma_0 \in T_{x_0} N \otimes T_{x_0} N$, we obtain

$$\beta(\varepsilon) = \frac{\varepsilon^2}{2}\{E[\Theta_0^{(2)}] + \bar{\Gamma}(y_\delta)\left(J\tau_0^\delta \Sigma_0 (J\tau_0^\delta)^T\right)\} + o(\varepsilon^2)$$

It follows that $\left.\dfrac{d\beta}{d\varepsilon}\right|_{\varepsilon = 0} = 0$, and hence that

$$\left.\frac{\partial}{\partial(\varepsilon^2)}\beta(\varepsilon)\right|_{\varepsilon = 0} = \frac{1}{2}\left.\frac{d^2\beta}{d\varepsilon^2}\right|_{\varepsilon = 0}$$

$$= \frac{1}{2}\{E[\Theta_0^{(2)}] + \bar{\Gamma}(y_\delta)\left(J\tau_0^\delta \Sigma_0 (J\tau_0^\delta)^T\right)\}$$

$$= \frac{1}{2}\{J\int_0^\delta \tau_s^\delta [D^2\xi(x_s)(\chi_s) - \Gamma(x_s)(\sigma \cdot \sigma(x_s))] ds + D^2\psi(x_\delta)(\chi_\delta) + \bar{\Gamma}(y_\delta)(J\chi_\delta J^T)$$

$$+ J\int_0^\delta \tau_s^\delta D^2\xi(x_s)(\Sigma_s) ds - J\tau_0^\delta \Gamma(x_0)(\Sigma_0) + D^2\psi(x_\delta)(\Sigma_\delta) + \bar{\Gamma}(y_\delta)(J\Sigma_\delta J^T) \}$$

where $\Sigma_s \equiv \tau_0^s \Sigma_0 (\tau_0^s)^T$. If we write

$$\Xi_s \equiv \chi_s + \Sigma_s, \tag{65}$$

then the formula becomes

$$\frac{1}{2}\{J\int_0^\delta \tau_s^\delta [D^2\xi(x_s)(\Xi_s)] ds + D^2\psi(x_\delta)(\Xi_\delta) + \bar{\Gamma}(y_\delta)(J\Xi_\delta J^T)\}$$

$$-\frac{1}{2}J\{\int_0^\delta \tau_s^\delta \Gamma(x_s)(\sigma \cdot \sigma(x_s)) ds + \tau_0^\delta \Gamma(x_0)(\Sigma_0)\} .$$

This establishes the formula (33). ◊

## 6.4 Step IV: Intrinsic Version of the Formula

It remains to prove that (24), with $\Pi_t$ as in (32), is the intrinsic version of (33). We abbreviate here by writing $\psi \bullet \phi_t$ as $\psi_t$. By definition of the flow of $\xi$,

$$\frac{\partial}{\partial t}\psi_t(x_0) = d\psi \bullet \xi(\phi_t(x_0)),$$

and so, differentiating with respect to $x$, and exchanging the order of differentiation,

$$\frac{\partial}{\partial t}D\psi_t(x_0) = D\left(\frac{\partial}{\partial t}\psi_t(x_0)\right) \tag{66}$$

or, by analogy with (29), taking $\theta_s^t \equiv D\psi_{t-s}(x_s) = D\psi(x_t) \bullet D(\phi_t \bullet \phi_s^{-1})(x_s) = D\psi(x_t)\tau_s^t$,

$$\frac{\partial \tau_s^t}{\partial t} = D\xi(x_t)\tau_s^t; \frac{\partial \theta_s^t}{\partial t} = D\psi(x_t) \bullet D\xi(x_t)\tau_s^t. \tag{67}$$

Since $\tau_t^\delta \tau_0^t = \tau_0^\delta$, we have $\theta_t^\delta \tau_0^t = D\psi(x_\delta)\tau_t^\delta \tau_0^t = \theta_0^\delta$, which upon differentiation yields

$$\left(\frac{\partial \theta_t^\delta}{\partial t}\right)\tau_0^t = -\theta_t^\delta\left(\frac{\partial \tau_0^t}{\partial t}\right) = -\theta_t^\delta D\xi(x_t)\tau_0^t,$$

which gives

$$\frac{\partial \tau_t^\delta}{\partial t} = -\tau_t^\delta D\xi(x_t); \frac{\partial \theta_t^\delta}{\partial t} = -\theta_t^\delta D\xi(x_t). \tag{68}$$

A further differentiation of (66) when $\psi$ is the identity yields

$$\frac{\partial}{\partial t}D^2\phi_t(x_0)\bigg)(v \otimes w) = D^2\xi(x_t)(\tau_0^t v \otimes \tau_0^t w) + D\xi(x_t)D^2\phi_t(x_0)(v \otimes w). \tag{69}$$

Combining (68) and (69), we have

$$\frac{\partial}{\partial t}(\theta_t^\delta D^2\phi_t(x_0)) = -\theta_t^\delta D\xi(x_t)D^2\phi_t(x_0) + \theta_t^\delta\{D^2\xi(x_t)(\tau_0^t(.) \otimes \tau_0^t(.)) + D\xi(x_t)D^2\phi_t(x_0)\}$$

$$\frac{\partial}{\partial t}(\theta_t^\delta D^2\phi_t(x_0))(v \otimes w) = \theta_t^\delta D^2\xi(x_t)(\tau_0^t v \otimes \tau_0^t w). \tag{70}$$

The formula (16) for $\nabla^\circ d\phi_t$ can be written as

$$\nabla^\circ d\phi_t(x_0)(v \otimes w) = D^2\phi_t(x_0)(v \otimes w) - \tau_0^t \Gamma(x_0)(v \otimes w) + \Gamma(x_t)(\tau_0^t v \otimes \tau_0^t w) \tag{71}$$

It is clear that

$$\frac{\partial}{\partial t}(\theta_t^\delta \tau_0^t \Gamma(x_0)) = \frac{\partial}{\partial t}(\theta_0^\delta \Gamma(x_0)) = 0.$$

Hence from (70) and (71) it follows that

$$\frac{\partial}{\partial t}\{(\psi_{\delta-t})_*\nabla^\circ d\phi_t(x_0)\}(v \otimes w) = \frac{\partial}{\partial t}\{\theta_t^\delta[D^2\phi_t(x_0)(v \otimes w) + \Gamma(x_t)(\tau_0^t v \otimes \tau_0^t w)]\}$$

$$= \theta_t^\delta D^2 \xi(x_t) (\tau_0^t v \otimes \tau_0^t w) + \frac{\partial}{\partial t} \{\theta_t^\delta \Gamma(x_t) (\tau_0^t v \otimes \tau_0^t w)\}. \tag{72}$$

The last term in (72) can be written, using (67), as

$$\frac{\partial}{\partial t} \{\theta_t^\delta \Gamma(x_t)\} (\tau_0^t v \otimes \tau_0^t w) + \theta_t^\delta \Gamma(x_t) \left(D\xi(x_t) \tau_0^t v \otimes \tau_0^t w + \tau_0^t v \otimes D\xi(x_t) \tau_0^t w\right), \tag{73}$$

for $v \otimes w \in T_{x_0} N \otimes T_{x_0} N$. We will replace $v \otimes w$ by

$$\Pi_t \equiv \tau_t^0 \Xi_t (\tau_t^0)^T = \Sigma_0 + \int_0^t \tau_s^0 (\sigma \cdot \sigma)(x_s) (\tau_s^0)^T ds \in T_{x_0} N \otimes T_{x_0} N \tag{74}$$

where $\Xi_t \equiv \chi_t + \tau_0^t \Sigma_0 (\tau_0^t)^T$. Observe that

$$\tau_0^t \Pi_t (\tau_0^t)^T = \Xi_t \in T_{x_t} N \otimes T_{x_t} N. \tag{75}$$

Moreover from (75) and (67), it is easily checked that

$$\frac{d\Xi_t}{dt} = (\sigma \cdot \sigma)(x_t) + D\xi(x_t) \Xi_t + \Xi_t \{D\xi(x_t)\}^T. \tag{76}$$

It follows from (72) - (76) that

$$\frac{\partial}{\partial t} \{(\psi_{\delta-t})_* \nabla^\circ d\phi_t(x_0) (\Pi_t)\} - (\psi_{\delta-t})_* \nabla^\circ d\phi_t(x_0) \left(\frac{d\Pi_t}{dt}\right) = \frac{\partial}{\partial t} \{(\psi_{\delta-t})_* \nabla^\circ d\phi_t(x_0)\} (\Pi_t)$$

$$= \left[\theta_t^\delta D^2 \xi(x_t) + \frac{\partial}{\partial t} \{\theta_t^\delta \Gamma(x_t)\}\right] (\Xi_t) + \theta_t^\delta \Gamma(x_t) [D\xi(x_t) \Xi_t + \Xi_t (D\xi(x_t))^T]$$

$$= \left[\theta_t^\delta D^2 \xi(x_t) + \frac{\partial}{\partial t} \{\theta_t^\delta \Gamma(x_t)\}\right] (\Xi_t) + \theta_t^\delta \Gamma(x_t) \left(\frac{d\Xi_t}{dt} - (\sigma \cdot \sigma)(x_t)\right)$$

$$= \theta_t^\delta D^2 \xi(x_t) (\Xi_t) - \theta_t^\delta \Gamma(x_t) ((\sigma \cdot \sigma)(x_t)) + \frac{\partial}{\partial t} \{\theta_t^\delta \Gamma(x_t) (\Xi_t)\}.$$

Since $\nabla^\circ d\phi_0 = 0$, it follows upon integration from 0 to $\delta$ that in $T_{\psi(x_\delta)} M$,

$$\psi_* \nabla^\circ d\phi_\delta(x_0) (\Pi_\delta) - \int_0^\delta (\psi_{\delta-t})_* (\nabla^\circ d\phi_t(x_0)) d\Pi_t =$$

$$D\psi(x_\delta) \{\int_0^\delta \tau_t^\delta [D^2\xi(x_t) (\Xi_t) - \Gamma(x_t)((\sigma \cdot \sigma)(x_t))] dt + \Gamma(x_\delta) (\Xi_\delta) - \tau_0^\delta \Gamma(x_0) (\Xi_0)\}. \tag{77}$$

However the formula (16) for $\nabla d\psi(x_\delta) (v \otimes w)$ can be written as

$$D^2 \psi(x_\delta) (v \otimes w) - J\Gamma(x_\delta) (v \otimes w) + \bar{\Gamma}(y_\delta) (Jv \otimes Jw). \tag{78}$$

where $y_\delta \equiv \psi(x_\delta)$, and $J \equiv D\psi(x_\delta)$. We take $v \otimes w = \Xi_\delta \in T_{x_\delta} N \otimes T_{x_\delta} N$, and add (77) and (78):

$$J\int_0^\delta \tau_t^\delta [D^2\xi(x_t)(\Xi_t) - \Gamma(x_t)((\sigma \cdot \sigma)(x_t))] \, dt + D^2\psi(x_\delta)(\Xi_\delta) - J\tau_0^\delta \Gamma(x_0)(\Xi_0) + \bar{\Gamma}(y_\delta)(J\Xi_\delta J^T)$$

$$= \nabla d\psi(x_\delta)((\phi_\delta)_* \Pi_\delta) + \psi_* \{\nabla^\circ d\phi_\delta(x_0)(\Pi_\delta) - \int_0^\delta (\phi_{\delta-t})_* (\nabla^\circ d\phi_t(x_0)) \, d\Pi_t\}.$$

The equivalence of (24) and (33) is established, completing the proofs of Theorems 4.5 and 4.7. ◊

## 7 The Canonical Sub-Riemannian Connection

The purpose of this section is to present a global geometric construction of a torsion-free connection $\nabla^\circ$ on the tangent bundle $TN$ which preserves, in some sense, a $C^2$ semi-definite metric $\langle . | . \rangle$ on the cotangent bundle $T^*N$ induced by a section $\sigma$ of $\text{Hom}(R^p; TN)$ of constant rank. In other words, we assume that there exists a rank $r$ vector bundle $E \to N$, a sub-bundle of the tangent bundle, such that $E_x = \text{range}(\sigma(x)) \subseteq T_x N$ for all $x \in N$.

Given such a section $\sigma$, we obtain a vector bundle morphism $\alpha : T^*N \to TN$ by the formula

$$\alpha(x) \equiv \sigma(x) \bullet \iota \bullet \sigma(x)^*, \quad x \in N, \tag{79}$$

where $\iota : (R^p)^* \to R^p$ is the canonical isomorphism induced by the Euclidean inner product. The relation between $\alpha$ and $\langle . | . \rangle$ is that, omitting $x$,

$$\mu \cdot \alpha(\lambda) \equiv \langle \mu | \lambda \rangle, \quad \forall \mu \in T^*N. \tag{80}$$

### 7.1 Lemma

*Under the constant-rank assumption, any Riemannian metric on $N$ induces an orthogonal splitting of the cotangent bundle of the form*

$$T_x^* N = \text{Ker}(\alpha(x)) \oplus F_x, \tag{81}$$

*where $F \to N$ is a rank $r$ sub-bundle of the cotangent bundle on which $\langle . | . \rangle$ is non-degenerate. There exists a vector bundle isomorphism $\alpha^\circ : T^*N \to TN$ such that $F_x = \text{Ker}(\alpha^\circ(x) - \alpha(x))$, and $\alpha^\circ(x)^{-1}$ is a generalized inverse to $\alpha(x)$, in the sense that*

$$\alpha(x) \bullet \alpha^\circ(x)^{-1} \bullet \alpha(x) = \alpha(x).$$

Proof: For any matrix $A$, $\text{range}(A) = \text{range}(AA^T)$, and so $\text{range}(\alpha(x)) = \text{range}(\sigma(x))$. It follows that

$$\dim \text{Ker}(\alpha(x)) = p - \dim E_x = p - r.$$

Given a Riemannian metric $g$ on $N$ (which always exists), let $\langle . | . \rangle^\circ$ be the dual metric on the cotangent bundle. Define

$$F_x \equiv \{\theta \in T_x^*N : \langle\theta|\lambda\rangle_x^\circ = 0 \ \forall \lambda \in \mathrm{Ker}\,(\alpha\,(x))\}\,.$$

We omit the proof that $F \to N$ is a vector bundle. Since $\dim \mathrm{Ker}\,(\alpha\,(x)) = p - r$, it follows that $\dim F_x = r$; since the rank of $\langle.|.\rangle$ is $r$, we see that $\langle\theta|\theta\rangle > 0$ for all non-zero $\theta \in F_x$. This shows that $\langle.|.\rangle$ is non-degenerate on the sub-bundle $F \to N$.

Now (81) results from the orthogonal decomposition of $T_x^*N$ with respect to $\langle.|.\rangle^\circ$. Hence an arbitrary $\lambda \in T_x^*N$ can be decomposed as $\lambda = \lambda_0 \oplus \lambda_1$, with $\lambda_0 \in \mathrm{Ker}\,(\alpha\,(x))$ and $\lambda_1 \in F_x$. The metric $\langle.|.\rangle^\circ$ induces a vector bundle isomorphism $\beta\colon T^*N \to TN$, namely

$$\mu \cdot \beta\,(\lambda) \equiv \langle\mu|\lambda\rangle^\circ,\ \forall \mu \in T^*N\,.$$

Now define $\alpha^\circ\colon T^*N \to TN$ by

$$\alpha^\circ\,(\lambda) \equiv \beta\,(\lambda_0) + \alpha\,(\lambda_1)\,.$$

It is clearly linear, and a vector bundle morphism. Since $\beta$ is injective,

$$\{\lambda \in T_x^*N : \alpha^\circ\,(\lambda) = \alpha\,(\lambda)\} = \{\lambda_0 \oplus \lambda_1 \in T_x^*N : \beta\,(\lambda_0) = 0\} = F_x,$$

which shows that $F_x = \mathrm{Ker}\,(\alpha^\circ\,(x) - \alpha\,(x))$. To show $\alpha^\circ$ is an isomorphism, it suffices to show that $\alpha^\circ\,(\lambda) \neq 0$ whenever $\lambda \neq 0$. When $\lambda_0 \neq 0$, non-degeneracy of $\langle.|.\rangle^\circ$ implies that

$$\lambda_0 \cdot \alpha^\circ\,(\lambda) = \langle\lambda_0|\lambda_0\rangle^\circ + \langle\lambda_0|\lambda_1\rangle = \langle\lambda_0|\lambda_0\rangle^\circ \neq 0\,.$$

On the other hand when $\lambda_0 = 0$, and $\lambda_1 \neq 0$, non-degeneracy of $\langle.|.\rangle$ on $F \to N$ implies that

$$\lambda_1 \cdot \alpha^\circ\,(\lambda) = \lambda_1 \cdot \alpha\,(\lambda_1) = \langle\lambda_1|\lambda_1\rangle \neq 0\,.$$

Hence $\alpha^\circ\colon T^*N \to TN$ is a vector bundle isomorphism as claimed. The generalized inverse property follows from the fact that $\alpha\,(x) \bullet \alpha^\circ\,(x)^{-1} \bullet \alpha\,(x)\,(\lambda) = \lambda_1$.    ◊

### 7.2    Proposition

*Suppose $\sigma$ is a constant-rank section of $\mathrm{Hom}\,(R^p;TN)$, inducing a semi-definite metric $\langle.|.\rangle$ on $T^*N$ and a vector bundle morphism $\alpha\colon T^*N \to TN$ as in (79) and (80). Suppose furthermore that $\alpha^\circ\colon T^*N \to TN$ is a vector bundle isomorphism such that $\alpha\,(x) \bullet \alpha^\circ\,(x)^{-1} \bullet \alpha\,(x) = \alpha\,(x)$, as in Lemma 7.1. Then $TN$ admits a canonical sub-Riemannian connection $\nabla^\circ$ for $\langle.|.\rangle$, with respect to $\alpha^\circ$, which is torsion-free, and such that the dual connection $\hat\nabla$ preserves $\langle.|.\rangle$ in the following sense: for vector fields $V$ in the range of $\alpha$, and for 1-forms $\theta, \lambda$ which lie in the sub-bundle $F \equiv \mathrm{Ker}\,(\alpha^\circ - \alpha)$,*

$$V\langle\theta|\lambda\rangle = \langle\hat\nabla_V\theta|\lambda\rangle + \langle\theta|\hat\nabla_V\lambda\rangle\,. \tag{82}$$

*[Here $\hat\nabla_Z\theta \cdot W = Z\,(\theta \cdot W) - \theta \cdot \nabla^\circ_Z W\,.$] For any 1-forms $\theta, \mu, \lambda$, and corresponding vector fields*

$$Y \equiv \alpha^\circ\,(\theta)\,,\ Z \equiv \alpha^\circ\,(\mu)\,,\ W \equiv \alpha^\circ\,(\lambda)\,,$$

*the formula for $\nabla°$ is:*

$$\mu \cdot \nabla°_Y W \equiv \frac{1}{2} \{ Y\langle\lambda|\mu\rangle + W\langle\mu|\theta\rangle - Z\langle\theta|\lambda\rangle + \lambda \cdot [Z, Y] + \mu \cdot [Y, W] - \theta \cdot [W, Z] \} . \tag{83}$$

### 7.2.a    Expression in Local Coördinates

Take local coördinates for $N$, so that $\alpha°(x)^{-1}$ is represented by a matrix $(g_{lm})$, and $\alpha(x)$ by a matrix $(\alpha^{jk})$, where by Lemma 7.2,

$$\sum_{k,r} \alpha^{jk} g_{kr} \alpha^{rm} = \alpha^{jm} .$$

Take $Y \equiv \partial/\partial x_i$, $W \equiv \partial/\partial x_j$, and $Z \equiv \partial/\partial x_k$ in (83), so that $\mu = \sum g_{ks} dx^s$, etc.; (83) becomes

$$\sum_s \Gamma^s_{ij} g_{sk} = \frac{1}{2} \sum_{r,s} \{ \frac{\partial}{\partial x_i}(g_{jr}\alpha^{rs} g_{sk}) + \frac{\partial}{\partial x_j}(g_{ir}\alpha^{rs} g_{sk}) - \frac{\partial}{\partial x_k}(g_{ir}\alpha^{rs} g_{sj}) \} . \tag{84}$$

When $\langle.|.\rangle$ is non-degenerate, then $\sum g_{jr}\alpha^{rs} = \delta^s_j$, and (84) reduces to the standard formula for the Levi-Civita connection for $g$, namely

$$\sum_s \Gamma^s_{ij} g_{sk} = \frac{1}{2} \left\{ \frac{\partial g_{jk}}{\partial x_i} + \frac{\partial g_{ik}}{\partial x_j} - \frac{\partial g_{ij}}{\partial x_k} \right\} .$$

### 7.2.b    Remark

A similar construction appears in formula (2.2) of Strichartz [27], where he cites unpublished work of N. C. Günther.

### 7.2.c    Proof of Proposition 7.2

First we check that the formula (83) defines a connection. The $R$-bilinearity of $(Y, W) \to \nabla°_Y W$ is immediate. To prove that $\nabla°_Y fW = f\nabla°_Y W + (Yf)W$ for all $f \in C^\infty(N)$, we replace $W$ by $fW$ and $\lambda$ by $f\lambda$ on the right side of (83), and the required identity holds. Verification that (83) is torsion-free is likewise a straightforward calculation.

Next we shall verify (82). By (80), and the definition of duality, taking $V \equiv \alpha(\theta)$,

$$\mu \cdot \nabla°_V \alpha(\lambda) = V(\mu \cdot \alpha(\lambda)) - \hat\nabla_V \mu \cdot \alpha(\lambda) = V\langle\lambda|\mu\rangle - \langle\lambda|\hat\nabla_V\mu\rangle .$$

Switch $\mu$ and $\lambda$ in the last expression to obtain:

$$\lambda \cdot \nabla°_V \alpha(\mu) = V\langle\lambda|\mu\rangle - \langle\mu|\hat\nabla_V\lambda\rangle .$$

Taking $U \equiv \alpha(\mu)$ and $S \equiv \alpha(\lambda)$, we see that

$$V\langle\lambda|\mu\rangle - \langle\lambda|\hat\nabla_V\mu\rangle - \langle\mu|\hat\nabla_V\lambda\rangle = \mu \cdot \nabla°_V S + \lambda \cdot \nabla°_V U - V\langle\lambda|\mu\rangle . \tag{85}$$

In terms of the splitting $T_x^* N = \text{Ker}(\alpha(x)) \oplus F_x$ of Lemma 7.1, we may write $\theta \equiv \theta_0 + \theta_1$, etc., and we find that, if $V \equiv \alpha(\theta)$, then $V = \alpha(\theta_1) = \alpha°(\theta_1)$, etc. It follows from (83) that

$$\mu_1 \cdot \nabla^\circ_V S \equiv \frac{1}{2} \{ V\langle\lambda_1|\mu_1\rangle + S\langle\mu_1|\theta_1\rangle - U\langle\theta_1|\lambda_1\rangle + \lambda_1 \cdot [U, V] + \mu_1 \cdot [V, S] - \theta_1 \cdot [S, U] \} ;$$

$$\lambda_1 \cdot \nabla^\circ_V U \equiv \frac{1}{2} \{ V\langle\lambda_1|\mu_1\rangle + U\langle\lambda_1|\theta_1\rangle - S\langle\theta_1|\mu_1\rangle + \mu_1 \cdot [S, V] + \lambda_1 \cdot [V, U] - \theta_1 \cdot [U, S] \} .$$

To prove (82), we can assume $\lambda_0 = \mu_0 = 0$, and now the right side of (85) becomes

$$\mu_1 \cdot \nabla^\circ_V S + \lambda_1 \cdot \nabla^\circ_V U - V\langle\lambda_1|\mu_1\rangle = 0,$$

as desired.                                                                                                         ◊

## 8    Future Directions

We would like to find out under what conditions on $\xi$, $\sigma \cdot \sigma$, $\psi$, and $h$ the system of PDE (17) - (18) has a unique solution for small time, other than the well-known case where $\sigma \cdot \sigma$ is non-degenerate, and $\xi = 0$; likewise for the parametrized family (22) - (23). It is likely that the conditions will involve the energy of the composite maps $\{\psi \bullet \phi_t, 0 \le t \le \delta\}$. Both stochastic and geometric methods should be considered. Another valuable project would be to derive bounds on the error of approximation involved in the linearization used in Theorem 4.7. This is likely to involve the curvature under the diffusion variance semi-definite metric – see Darling [9].

**Acknowledgments:** The author thanks the Statistics Department at the University of California at Berkeley for its hospitality during the writing of this article, Anton Thalmaier for allowing access to unpublished work, and James Cloutier, James Eells, Etienne Pardoux, and Richard Schoen for their advice.